\documentclass[12pt]{article}
\usepackage{amssymb}
\usepackage{amsmath}

\begin{document}

\begin{center}
\textbf{\large Homeomorphic Changes of Variable\\ and Fourier
Multipliers}
\end{center}

\begin{center}
\textsc{Vladimir Lebedev and Alexander Olevskii}
\end{center}

\begin{quotation}
{\small \textsc{Abstract.} We consider the algebras $M_p$ of
Fourier multipliers and show that for every bounded continuous
function $f$ on $\mathbb R^d$ there exists a self-homeomorphism
$h$ of $\mathbb R^d$ such that the superposition $f\circ h$ is
in $M_p(\mathbb R^d)$ for all $p$, $1<p<\infty$. Moreover,
under certain assumptions on a family $K$ of continuous
functions, one $h$ will suffice for all $f\in K$. A similar
result holds for functions on the torus~$\mathbb T^d$. This may
be contrasted with the known solution of Luzin's problem
related to the Wiener algebra.


}
\end{quotation}

\quad

\begin{center}
\textbf{1. Introduction}
\end{center}

To what extent the behavior of the Fourier series of a
continuous function can be improved by a change of variable?
The problem originates from a theorem of Bohr and P\'al (see
[1, Ch.~4, Sec.~12], [9], [23]), who showed that given a
continuous real-valued function $f$ on the circle $\mathbb
T=\mathbb R/2\pi \mathbb Z$, there exists a self-homeomorphism
$h$ of $\mathbb T$ such that the Fourier series of the
superposition $f\circ h$ converges uniformly. In addition, the
proof yields a condition on the decay of the Fourier
coefficients of $f\circ h$; namely,
$$
\widehat{f\circ h}\in l^p(\mathbb Z)
\eqno(1)
$$
for all $p>1$.

This result inspired the following problem posed by N. Lusin
(see [1, Ch.~4, Sec.~12], [8, Ch. VII, Sec. 9]): Is it possible
to attain condition (1) for $p=1$? In other words, is it true
that for every continuous function on $\mathbb T$ there exists
a change of variable which brings it into the Wiener
algebra~$A(\mathbb T)$?

It is worth noting that the original proof of the Bohr--P\'al
theorem involves Riemann's theorem on conformal mappings. A
purely real-analytic proof was found in [25].

A solution (in the negative) of Lusin's problem was obtained in
[22], where it is shown that there exists a continuous
real-valued function $f$ on $\mathbb T$ such that $f\circ h
\notin A(\mathbb T)$ whenever $h$ is a self-homeomorphism of
$\mathbb T$. Some ideas of P. Cohen were used in the proof.
Simultaneously, the same result for a complex-valued $f$ was
obtained in [10] (this weaker result amounts to the fact that,
in general, there is no single change of variable which will
bring two continuous real-valued functions into~$A(\mathbb
T)$).

Subsequently, for certain function spaces, naturally arising in
harmonic analysis, the question of whether every continuous
function can be transformed by a suitable homeomorphic change
of variable into a function that belongs to a given space, was
studied by various authors. Some of these studies concern the
possibility of simultaneous improvement of several functions by
means of a single change of variable; we note in particular the
work [11] where the Bohr--P\'al theorem is extended to compact
families of continuous functions.

For a survey on the subject see [23], [9]. For more recent
results see [24], [4, Ch. 9], [15], [16], [18], [19]. We also
mention the papers [14], [13] which study the growth of the
partial sums of the Fourier series of $f\circ h$ for random
homeomorphisms $h$.

In the present paper we investigate the problem on changes of
variable in relation with the Fourier Multiplier Algebras. We
consider multipliers on $\mathbb R^d$ and on the torus $\mathbb
T^d$.

Let $G$ be one of the groups $\mathbb R^d$ or $\mathbb Z^d$ and
let $\Gamma$ be its dual group, i.e., $\Gamma=\mathbb R^d$ or
$\Gamma=\mathbb T^d$, correspondingly. Consider a function
$m\in L^\infty (\Gamma)$ and the operator $Q$ defined by
$$
\widehat{Qf}=m\cdot\widehat{f}, \qquad f\in L^p\cap L^2(G),
\eqno(2)
$$
where $\widehat{~~}$ stands for the Fourier transform on~$G$.
The function $m$ is called an $L^p$ multiplier ($1\leq
p\leq\infty$) if
$$
\|Qf\|_{L^p(G)}\leq c \|f\|_{L^p(G)}, \qquad f\in L^p\cap L^2(G),
\eqno(3)
$$
where $c>0$ does not depend on $f$. The space of all such
multipliers is denoted by $M_p(\Gamma)$. The norm on
$M_p(\Gamma)$ is defined by setting $\|m\|_{M_p(\Gamma)}$ to be
the smallest $c$ for which (3) holds. The space $M_p(\Gamma)$
equipped with this norm is a Banach algebra with the usual
multiplication of functions. Clearly, if $p<\infty$, then the
operator $Q$ that corresponds to the function $m$ can be
uniquely extended to a bounded operator on $L^p(G)$, and
retaining the notation $Q$ for this extension, we have
$\|Q\|_{L^p(G)\rightarrow L^p(G)}=\|m\|_{M_p(\Gamma)}$. We also
note that the operator $Q$ is translation-invariant. The
converse also holds for $p<\infty$: every translation-invariant
bounded operator on $L^p(G)$ has the form (2), where $m\in
M_p(\Gamma)$. For basic properties of multipliers see [3], [5].

It is known that $M_2$ coincides with $L^\infty$. It is also
known that $M_1=M_\infty$, and, at the same time, $M_1(\mathbb
T^d)$ coincides with the Wiener algebra $A(\mathbb T^d)$ and
$M_1(\mathbb R^d)$ coincides with the algebra $B(\mathbb R^d)$
of the Fourier transforms of (complex) bounded regular Borel
measures on $\mathbb R^d$. Note that the negative solution of
Luzin's problem immediately implies a similar result for
functions on the real line: there exists a bounded continuous
real-valued function $f$ on $\mathbb R$ such that $f\circ
h\notin B(\mathbb R)$ for every self-homeomorphism $h$ of
$\mathbb R$ (for details, see Sec. 5). Thus, in general, there
is no change of variable which will bring a continuous
real-valued function on $\mathbb T$ (a bounded continuous
function on $\mathbb R$) into $M_1=M_\infty$.

In this paper we show that for every bounded continuous
function $f$ on $\mathbb R^d$ (for every continuous function
$f$ on $\mathbb T^d$) there is a homeomorphic change of
variable $h$ such that $f\circ h\in\bigcap_{1<p<\infty}M_p$.
Moreover, under certain natural assumptions on a function
family $K$, one change of variable will suffice for all $f\in
K$. An important role in the proof is played by a result of
Sj\"ogren and Sj\"olin on Littlewood--Paley partitions.

The exact statements of our results are given in Section 2
below. Section 3 contains preliminaries. The proofs are given
in Section 4. The concluding Section 5 contains several remarks
and open problems; in particular, we focus on the famous
Beurling--Helson theorem and briefly discuss some recent
results and open questions related to this theorem and to its
version for Multiplier Algebras.

\quad

\begin{center}
\textbf{2. Statement of Results}
\end{center}

Let $f$ be a function defined on a set $E$ (we assume
that $E\subseteq\mathbb R^d$ or $E\subseteq\mathbb T^d$). By
$\omega(f, E, \delta)$ we denote the modulus of continuity of
$f$ on $E$:
$$
\omega(f, E, \delta)=\sup_{\underset{|t_1-t_2|\leq
\delta}{t_1, t_2\in E}} |f(t_1)-f(t_2)|,
\qquad \delta\geq 0
$$
($|x|$ stands for the length of a vector $x\in\mathbb R^d$).
Let $K$ be a family of functions on $E$. We define the modulus
of continuity $\omega(K, E, \delta)$ of $K$ on $E$ by
$$
\omega(K, E, \delta)=
\sup_{f\in K} \, \omega(f, E, \delta), \qquad\delta\geq 0.
$$
We say that $K$ is uniformly equicontinuous on $E$ if
$\omega(K, E, \delta)\rightarrow 0$ as $\delta\rightarrow+0$.

By $C(\mathbb R^d)$ and $C(\mathbb T^d)$ we denote the classes
of continuous functions on $\mathbb R^d$ and $\mathbb T^d$,
respectively. The results of this paper are the following two
theorems.

\quad

\textsc{Theorem 1.} \emph{Let $f\in C(\mathbb R^d)$ be a
bounded function. Then there exists a self-homeomorphism $h$ of
$\mathbb R^d$ such that $f\circ h\in
\bigcap_{1<p<\infty}M_p(\mathbb R^d)$. Moreover, if a family
$K\subseteq C(\mathbb R^d)$ of bounded functions is uniformly
equicontinuous on every ball in $\mathbb R^d$, then there
exists a self-homeomorphism $h$ of $\mathbb R^d$ such that
$f\circ h\in \bigcap_{1<p<\infty}M_p(\mathbb R^d)$ for all
$f\in K$}.

\quad

\textsc{Theorem 2.} \emph{Let $f\in C(\mathbb T^d)$. Then there
exists a self-homeomorphism $h$ of $\mathbb T^d$ such that
$f\circ h\in \bigcap_{1<p<\infty}M_p(\mathbb T^d)$. Moreover,
if a family $K\subseteq C(\mathbb T^d)$ is uniformly
equicontinuous on $\mathbb T^d$, then there exists a
self-homeomorphism $h$ of $\mathbb T^d$ such that $f\circ h\in
\bigcap_{1<p<\infty}M_p(\mathbb T^d)$ for all $f\in K$.}

\quad

\begin{center}
\textbf{3. Preliminaries}
\end{center}

In this section we recall some notions and facts from the
Littlewood--Paley and multiplier theory.

\quad

A. For $1/p+1/q=1$ we have $M_p(\Gamma) = M_q(\Gamma)$ and
$\|\cdot\|_{M_p}=\|\cdot\|_{M_q}$. If $1\leq p_1<p_2\leq 2$,
then $M_{p_1}(\Gamma)\subset M_{p_2}(\Gamma)$ and
$\|\cdot\|_{M_{p_2}}\leq \|\cdot\|_{M_{p_1}}$ (see, e.g., [3]).

\quad

B. The indicator function $1_I$ of any rectangle
$I\subseteq\mathbb R^d$ belongs to $M_p(\mathbb R^d)$ for all
$p, 1<p<\infty$, see, e.g., [3]. (As usual, we set $1_I(t)=1$
if $t\in I$ and $1_I(t)=0$ if $t\notin I$.)

\quad

C. Let $m\in M_p(\mathbb R^d)$, $1\leq p\leq\infty$, and let $l
\colon \mathbb R^d\rightarrow\mathbb R^d$ be a nondegenerate
affine mapping; then $m\circ l\in M_p(\mathbb R^d)$ and
$\|m\circ l\|_{M_p(\mathbb R^d)}=\|m\|_{M_p(\mathbb R^d)}$
(see, e.g., [5, Ch.~I, Sec.~1.3]).

\quad

D. Let $m_0$ be a function in $M_p(\mathbb R^d)$, $1\leq
p\leq\infty$, vanishing outside a cube $I$ with edges of length
$2\pi$ parallel to the coordinate axes. Consider the extension
$m$ of $m_0$ to $\mathbb R^d$ which is $2\pi$-periodic in each
variable. The function $m$ belongs to $M_p(\mathbb T^d)$, and $
\|m\|_{M_p(\mathbb T^d)}\leq c_p\|m_0\|_{M_p(\mathbb R^d)}, $
where $c_p$ is a positive constant not depending on $m_0$. This
theorem on extension of multipliers is well known [6,
Theorem~2.3] (for $p=1, \infty$ the result follows from the
local properties of the Fourier transforms of measures, see.,
e.g., [8, Ch. II, Sec.~4]).

\quad

E. Let $I\subseteq\mathbb R^d$ be a rectangle. By $S_I$ we
denote the operator on $L^p(\mathbb R^d)$ that corresponds to
multiplication by $1_I$, i.e., the operator defined by
$$
\widehat{S_I(f)}=1_I\cdot\widehat{f}, \qquad f\in L^p\cap L^2(\mathbb R^d).
$$
Let $\Delta$ be a family of rectangles which form a partition
of $\mathbb R^d$, i.e., a family of pairwise disjoint
rectangles in $\mathbb R^d$ such that the complement $\mathbb
R^d\setminus\bigcup_{I\in\Delta} I$ has Lebesgue measure zero.
Consider the corresponding Littlewood--Paley square function
$S^\Delta(f)$:
$$
S^\Delta(f)=\bigg(\sum_{I\in\Delta} |S_I(f)|^2 \bigg)^{1/2}.
$$
The partition $\Delta$ is called an $\mathrm{LP}$ partition if,
for all $p$, $1<p<\infty$, we have
$$
a_p\cdot \|f\|_{L^p(R^d)}\leq \|S^\Delta(f)\|_{L^p(R^d)}\leq b_p\cdot
\|f\|_{L^p(R^d)},
$$
where $a_p=a_p(\Delta)$ and $b_p=b_p(\Delta)$ are positive
constants independent of $f$. A classical example of an
$\mathrm{LP}$ partition of the real line $\mathbb R$ is the
dyadic partition, that is the family of the intervals $I_k,
\,k\in\mathbb Z,$ of the form $I_k=(2^{k-1}, 2^k)$ for $k=1, 2,
\ldots$, $I_0=(-1, 1)$, and $I_k=(-2^{-k}, -2^{-k-1})$ for
$k=-1, -2, \ldots$.

Let $\Delta$ be a family of rectangles which form an
$\mathrm{LP}$ partition of $\mathbb R^d$. Consider an arbitrary
function $m\in L^\infty(\mathbb R^d)$ constant on each
rectangle $I\in\Delta$. Then $m$ is in $M_p(\mathbb R^d)$ for
all $p$, $1<p<\infty$, and
$$
\|m\|_{M_p(\mathbb R^d)}\leq c(p, \Delta)\cdot\|m\|_{L^\infty(\mathbb R^d)},
$$
where $c(p, \Delta)>0$ does not depend on $m$ (see, e.g., [3,
1.2.7, 1.2.8]).

Given a  family $\Delta$ of intervals in $\mathbb R$, let
$\Delta_d$ denote the family of rectangles in $\mathbb R^d$
generated by $\Delta$, that is, the family of all rectangles
$I$ of the form $I=I_1\times I_2\times \ldots I_d$, where each
factor belongs to $\Delta$. If $\Delta$ is an $\mathrm{LP}$
partition of $\mathbb R$, then $\Delta_d$ is an $\mathrm{LP}$
partition of $\mathbb R^d$ (see [3, Theorem~1.3.4]).

\quad

F. By the dyadic partition of the interval $(0, 1)$ we mean the
family of the intervals $I_k, \,k\in\mathbb Z$, defined by
$$
I_k=(1-2^{-k-1}, \,1-2^{-k-2}), \quad k=1,2,\ldots;
$$
$$
I_0=(1/4, \,3/4);
$$
$$
I_k=(2^{k-2}, \,2^{k-1}), \quad k=-1,-2,\ldots\,.
$$
We extend this definition to any interval $(a, b)$ in the
obvious way by translation and rescaling. Suppose now that
$\Delta$ is a family of intervals which form an $\mathrm{LP}$
partition of $\mathbb R$. Then the family of all intervals,
obtained by dyadic partition of each interval $I\in\Delta$ is
an $\mathrm{LP}$ partition of $\mathbb R$ as well. This
immediately follows from a result of Sj\"ogren and Sj\"olin
[26, Theorem~1.2].

\quad

\begin{center}
\textbf{4. Proofs of the Theorems}
\end{center}

First, we prove a simple auxiliary lemma.

\quad

\textsc{Lemma.} \emph{Let $f\in C(\mathbb R^d)$. Then there
exists a self-homeomorphism $\psi$ of $\mathbb R^d$ such that
$f\circ\psi$ is uniformly continuous on $\mathbb R^d$.
Moreover, if a family $K\subseteq C(\mathbb R^d)$ is uniformly
equicontinuous on every ball in $\mathbb R^d$, then there
exists a self-homeomorphism $\psi$ of $\mathbb R^d$ such that
the family $\{f\circ \psi, f\in K\}$ is uniformly
equicontinuous on $\mathbb R^d$.}

\quad

\emph{Proof.} We give the proof in the general case of families
of functions. Consider spherical shells
$$
\Theta_j=\{x\in\mathbb R^d : j\leq |x|\leq j+1\}, \quad j=0, 1, 2,
\ldots
$$
($\Theta_0$ is a ball). Let $\omega_j$ be the modulus of
continuity of the family $K$ on $\Theta_j$, i.e.,
$$
\omega_j(\delta)=\sup_{f\in K}\omega (f, \Theta_j, \delta).
$$
By assumption, $\omega_j(\delta)\rightarrow 0$ as
$\delta\rightarrow +0$ for each $j$. So, we can find a
decreasing sequence $b_j, \,j=0, 1, 2, \ldots,$ such that
$0<b_j<1$ and
$$
\omega_j(b_j)\rightarrow 0 \quad \textrm{as} \quad j\rightarrow\infty.
\eqno(4)
$$

We set
$$
a_j=\frac{b_j}{j+2}, \qquad j=1, 2, \ldots.
$$
The sequence $a_j, \,j=1, 2, \ldots,$ decreases, and $0<a_j<1$.
Define numbers $r_j$, $j=0, 1, 2, \ldots\,$,  by
$$
r_0=0; \qquad r_j=\sum_{s=1}^j \frac{1}{a_s}, \quad j=1, 2, \ldots\,.
\eqno(5)
$$
Clearly, there exists a function $g$ on $[0, +\infty)$ with
the following properties:

(i) $g(0)=0$;

(ii)  $g$ is continuous and strictly increasing;

(iii) for each $j=0, 1, 2, \ldots$ the function $g$ is linear
on the interval $[r_j, r_{j+1}]$ and maps $[r_j, r_{j+1}]$ onto
the interval $[j, j+1]$;

(iv) the slope of $g$ on $[r_j, r_{j+1}]$ equals $a_{j+1}$ (see
(5)).

For  $x\in\mathbb R^d$, we set
$$
\psi(x)=g(|x|)\frac{x}{|x|}, \quad x\neq 0; \qquad \psi(0)=0.
$$
One can easily see that $\psi$ is a self-homeomorphism of
$\mathbb R^d$.

Consider the spherical shells
$$
\Omega_j=\{x : r_j\leq |x|\leq
r_{j+1}\}, \quad j=0, 1, 2, \ldots
$$
($\Omega_0$ is a ball). Note, that the image of $\Omega_j$
under $\psi$ is the shell $\Theta_j$. We claim that
$$
|\psi(x)-\psi(y)|\leq b_j|x-y| \quad \textrm{for all}\quad x, y\in \Omega_j,
\qquad j=0, 1, 2, \ldots.
\eqno(6)
$$
Since $\psi(x)=a_1 x$ on $\Omega_0$, we see that if $x, y\in
\Omega_0$, then
$$
|\psi(x)-\psi(y)|=a_1|x-y|=(b_1/3)|x-y|\leq b_0|x-y|.
$$
Let $j\geq 1$. Observe that the mapping $\gamma(x)=x/|x|$ has
the property that $|\gamma(x_1)-\gamma(x_2)|\leq(1/r)|x_1-x_2|$
for all $x_1, x_2$ that lie outside the ball $B_r=\{x :
|x|<r\}, \,r>0$. So, if $x, y\in \Omega_j$, then
$$
|\gamma(x)-\gamma(y)|\leq\frac{1}{r_j}|x-y|.
$$
At the same time
$$
|g(|x|)-g(|y|)|=a_{j+1}||x|-|y||\leq a_{j+1}|x-y|.
$$
Whence, taking into account that $g(|y|)\leq j+1$ and $r_j\geq
1/a_j$, we obtain
$$
|\psi(x)-\psi(y)|=
\big|\big(g(|x|)-g(|y|)\big)\gamma(x)+g(|y|)
\big(\gamma(x)-\gamma(y)\big)\big|
$$
$$
\leq a_{j+1}|x-y|+(j+1)\frac{1}{r_j}|x-y|
\leq (j+2)a_j|x-y|=b_j|x-y|.
$$
Thus, (6) holds.

Let $f\in K$. We shall estimate the modulus of continuity of
the superposition $f\circ\psi$ on $\mathbb R^d$. Let
$0<\delta\leq 1$, and let $x, y\in\mathbb R^d, \,|x-y|\leq
\delta$. Since the thickness of each shell $\Omega_j$ is
greater than $1$ (it equals $1/a_{j+1}$), we see that either
the points $x$ and $y$ belong to the same shell, say
$\Omega_j$, or there are two neighboring shells $\Omega_j$ and
$\Omega_{j+1}$ such that $x$ is in one of them and $y$ is in
the other one. Consider the first case when $x$ and $y$ are in
$\Omega_j$. Then the points $\psi(x), \psi(y)$ belong to the
shell $\Theta_j$. So (see (6)),
$$
|f\circ\psi(x)-f\circ\psi(y)|\leq
\omega_j(|\psi(x)-\psi(y)|)\leq \omega_j(b_j\delta).
$$
Consider the second case when $x\in \Omega_j$ and $y\in
\Omega_{j+1}$. The line segment that joins $x$ and $y$ contains
a point $z$ that belongs to both $\Omega_j$ and $\Omega_{j+1}$.
Using the estimate obtained in the first case, we have
$$
|f\circ\psi(x)-f\circ\psi(y)|\leq|f\circ\psi(x)-f\circ\psi(z)|+
|f\circ\psi(z)-f\circ\psi(y)|\leq\omega_j(b_j\delta)+
\omega_{j+1}(b_{j+1}\delta).
$$
Thus, we see that for $0<\delta\leq 1$
$$
\omega(f\circ\psi, \mathbb R^d, \delta)\leq 2\sup_{j\geq 0}\omega_j(b_j\delta).
$$

To complete the proof it remains to note that (4) implies
$$
\sup_{j\geq 0}\omega_j(b_j\delta)\rightarrow 0 \quad\textrm{as}\quad\delta\rightarrow+0.
$$

\quad

\emph{Proof of Theorem 1.} Let $I=(l, r)$ be a bounded or
unbounded interval in $\mathbb R$. We say that intervals
$I_k\subset I$, $k\in\mathbb Z$, form an ordered partition of
$I$ if $I_k=(\theta_k, \theta_{k+1})$, where
$\theta_k<\theta_{k+1}$ for all $k\in\mathbb Z$,
$\lim_{k\rightarrow+\infty}\theta_k=r$, and
$\lim_{k\rightarrow-\infty}\theta_k=l$.

Suppose that intervals $I_{s_1}$, $s_1\in\mathbb Z$, form an
ordered partition of $\mathbb R$. For each fixed $s_1$,  let
$I_{s_1, s_2}$, $s_2\in\mathbb Z,$ be certain intervals, which
form an ordered partition of $I_{s_1}$, and for each $\nu$ and
integers $s_1, s_2, \ldots, s_\nu$ let $I_{s_1, s_2, \ldots,
s_\nu, s_{\nu+1}}$, $s_{\nu+1}\in\mathbb Z$, be intervals,
which form an ordered partition of $I_{s_1, s_2, \ldots,
s_\nu}$. Proceeding, we obtain a certain family of intervals:
$$
\{I_{s_1, s_2, \ldots, s_\nu} : \nu=1, 2, \ldots,
\quad\textrm{and}\quad s_1, s_2, \ldots, s_\nu\in\mathbb Z\}.
$$
We refer to any family of intervals thus obtained as a net. For
each fixed $\nu$, the intervals $I_{s_1, s_2, \ldots, s_\nu}$
are called intervals of rank $\nu$ (of a given net).

Given a set $E$ and a function $f$ on $E$, by $\mathrm{osc}_E
f$ we denote the oscillation of $f$ on $E$:
$\mathrm{osc}_E\,f=\sup_{t_1, t_2\in E}|f(t_1)- f(t_2)|$.

We shall construct two nets of intervals. Clearly, the dyadic
partition of the real line or of an interval (see Section~3,\,E
and F) is an ordered partition. Consider the intervals
$I_{s_1}$, $s_1\in\mathbb Z$, which form the dyadic partition
of $\mathbb R$. If intervals $I_{s_1, s_2, \ldots, s_\nu}$ of
rank $\nu$ are already defined, then we define $I_{s_1, s_2,
\ldots, s_\nu, s_{\nu+1}}$, $s_{\nu+1}\in\mathbb Z$, to be the
intervals that form the dyadic partition of $I_{s_1, s_2,
\ldots, s_\nu}$. Thereby, we have constructed the first net,
which we denote by $\alpha$.

By $\alpha_{d}(\nu)$ we denote the family of all rectangles in
$\mathbb R^d$ obtained as the Cartesian product of any $d$
intervals of rank $\nu$ of the net $\alpha$. It follows from
the properties of $\mathrm{LP}$ partitions listed in
Section~3,\,E and F, that if $m$ is a function in
$L^\infty(\mathbb R^d)$ constant on each rectangle which
belongs to $\alpha_d(\nu)$, then $m\in M_p(\mathbb R^d)$ for
all $p$, $1<p<\infty$, and
$$
\|m\|_{M_p(\mathbb R^d)}\leq c(p, \nu)\cdot
\|m\|_{L^\infty(\mathbb R^d)}, \qquad 1<p<\infty,
$$
where $c(p, \nu)>0$ may depend only on $p$ and $\nu$. In what
follows, we assume that
$$
c(p, \nu)=\sup_{m\in\mathcal{P}_\nu, \, m\neq 0}
\frac{\|m\|_{M_p(\mathbb R^d)}}{\|m\|_{L^\infty(\mathbb R^d)}},
\eqno(7)
$$
where $\mathcal{P}_\nu$ is the class of all functions in
$L^\infty(\mathbb R^d)$ constant on each rectangle from
$\alpha_d(\nu)$.

We proceed to the construction of the second net. By using the
Lemma, we can assume that $K$ is a family of bounded functions
uniformly equicontinuous on the whole $\mathbb R^d$. Let
$\omega$ be the modulus of continuity of $K$. We have
$$
\sup_{\underset{|t_1-t_2|\leq
\delta}{t_1, t_2\in\mathbb R^d}}|f(t_1)-f(t_2)|\leq
\omega(\delta), \quad \delta>0, \qquad\textrm{for all} \,\, f\in K,
$$
where $\omega(\delta)\rightarrow 0$ as $\delta\rightarrow +0$,
and $\omega$ is nondecreasing on $[0, +\infty)$. Fix a positive
decreasing sequence $\delta_\nu$, $\nu=1, 2, \ldots\,$, that
tends to $0$ so fast that
$$
\sum_{\nu=2}^\infty c\bigg(1+\frac{1}{\nu},
\,\nu\bigg)\omega(\delta_{\nu-1}\sqrt{d})<\infty.
\eqno(8)
$$
Note that (8) implies that
$$
\sum_{\nu=2}^\infty c(p, \nu)\omega(\delta_{\nu-1}\sqrt{d})<\infty
\eqno(9)
$$
for all $p$, $1<p<\infty$. Indeed, it suffices to observe that
if $\nu$ is large enough, then $c(p, \nu)\leq c(1+1/\nu, \nu)$
(see Section~3,\,A).

Let $J_{s_1}$, $s_1\in\mathbb Z$, be intervals of length at
most $\delta_1$ which form an ordered partition of the line
$\mathbb R$. For a fixed $\nu$, assuming that all the intervals
$J_{s_1, s_2, \ldots, s_\nu}$ are already defined, consider an
ordered partition of each interval $J_{s_1, s_2, \ldots,
s_\nu}$ by intervals $J_{s_1, s_2, \ldots, s_\nu, s_{\nu+1}}$,
$s_{\nu+1}\in\mathbb Z$, of length at most $\delta_{\nu +1}$.
Thereby, we obtain the second net, which we denote by $\beta$.

Clearly, there exists a self-homeomorphism  $\varphi$ of
$\mathbb R$ such that $\varphi(I_{s_1, s_2, \ldots,
s_\nu})=J_{s_1, s_2, \ldots, s_\nu}$ for all $\nu=1,2, \ldots$
and $s_1, s_2, \ldots, s_\nu$ (the intervals $I_{s_1, s_2,
\ldots, s_\nu}$ and $J_{s_1, s_2, \ldots, s_\nu}$ belong to the
nets $\alpha$ and $\beta$, respectively). We define a
homeomorphism $h$ of $\mathbb R^d$ onto itself by
$$
h(t)=(\varphi(t_1), \varphi(t_2), \ldots, \varphi(t_d)),
\qquad t=(t_1, t_2, \ldots, t_d)\in\mathbb R^d.
\eqno(10)
$$

Consider an arbitrary function $f\in K$ and set $g=f\circ h$.
Let us verify that $g\in\bigcap_{1<p<\infty}M_p(\mathbb R^d)$.
Given a rectangle $I\subset\mathbb R^d$ denote the center of
$I$ by $c_I$. For $\nu=1, 2, \ldots$, let $g_\nu$ be the
function that takes the constant value $g(c_I)$ on each
rectangle $I\in\alpha_d(\nu)$. Clearly,
$$
\|g-g_\nu\|_{L^\infty(\mathbb R^d)}\leq\sup_{I\in\alpha_d(\nu)}\mathrm{osc}_I\,g.
\eqno(11)
$$
Note that if $I\in\alpha_d(\nu)$, then the image $h(I)$ of $I$
under $h$ is a certain rectangle whose edges are of length at
most $\delta_\nu$, which implies that $\mathrm{diam}\,h(I)\leq
\delta_\nu\sqrt{d}$. So,
$$
\mathrm{osc}_I\,g=\mathrm{osc}_{h(I)}\,f\leq
\omega(f, \delta_\nu\sqrt{d})\leq
\omega(\delta_\nu\sqrt{d}) \quad \textrm{for all}\,\,
I\in\alpha_d(\nu),\quad \nu=1, 2, \ldots.
$$
Therefore (see (11)),
$$
\|g-g_\nu\|_{L^\infty(\mathbb R^d)}\leq\omega(\delta_\nu\sqrt{d}).
\eqno(12)
$$
Hence, for $\nu\geq 2$, we obtain
$$
\|g_{\nu}-g_{\nu-1}\|_{L^\infty(\mathbb R^d)}\leq
\|g_{\nu}-g\|_{L^\infty(\mathbb R^d)}+\|g-g_{\nu-1}\|_{L^\infty(\mathbb R^d)}
\leq 2\omega(\delta_{\nu-1}\sqrt{d}).
\eqno(13)
$$
Let $1<p<\infty$. Since $g_\nu\in\mathcal{P_\nu}$, it follows
that $g_\nu\in M_p(\mathbb R^d)$, $\nu=1, 2, \ldots\,$. Note
that $g_{\nu}-g_{\nu-1}\in\mathcal{P}_{\nu}$ for $\nu\geq 2$;
therefore (see (7) and (13)), we have
$$
\|g_{\nu}-g_{\nu-1}\|_{M_p(\mathbb R^d)}\leq c(p,
\nu)2\omega(\delta_{\nu-1}\sqrt{d}).
$$
Thus,
$$
\|g_{n+m}-g_n\|_{M_p(\mathbb R^d)}\leq
\sum_{\nu=n+1}^{n+m}\|g_{\nu}-g_{\nu-1}\|_{M_p(\mathbb R^d)}
\leq\sum_{\nu=n+1}^{n+m} 2c(p,
\nu)\omega(\delta_{\nu-1}\sqrt{d}).
\eqno(14)
$$
Taking (9) into account, we see that the sequence $g_\nu$,
$\nu=1, 2, \ldots\,$, converges in $M_p(\mathbb R^d)$ (recall
that $M_p$ is a Banach space). At the same time from (12) it
follows that this sequence converges to $g$ in
$L^\infty(\mathbb R^d)$. It remains to recall that
$\|\cdot\|_{L^\infty}=\|\cdot\|_{M_2}\leq\|\cdot\|_{M_p}$. This
completes the proof of Theorem 1.

\quad

\emph{Proof of Theorem 2.} Let $K$ be a uniformly
equicontinuous family of functions on $\mathbb T^d$. We
identify each function on $\mathbb T^d$ with a  $2\pi$-periodic
(in each variable) function on $\mathbb R^d$ in the standard
manner. Thus, $K$ is a family of bounded functions uniformly
equicontinuous on $\mathbb R^d$. Following the proof of
Theorem~1 (now we do not need the Lemma), we obtain a
self-homeomorphism $\varphi$ of $\mathbb R$ such that defining
self-homeomorphism $h$ of $\mathbb R^d$ by (10) we have $f\circ
h\in\bigcap_{1<p<\infty}M_p(\mathbb R^d)$ for all $f\in K$.
Consider the interval $J=\varphi^{-1}([0, 2\pi])$ which is the
preimage of $[0, 2\pi]$ under $\varphi$. Let $l$ be an affine
self-mapping of the real line for which $l([0, 2\pi])=J$. We
set $\varphi_1=\varphi\circ l$ and
$$
h_1(x)=(\varphi_1(x_1), \varphi_1(x_2), \ldots, \varphi_1(x_d)),
\quad x=(x_1, x_2, \ldots, x_d)\in\mathbb R^d.
$$
Using the assertion on superpositions of multipliers with
affine mappings (see Section~3,\,C), we conclude that $f\circ
h_1\in\bigcap_{1<p<\infty}M_p(\mathbb R^d)$ for all $f\in K$.
Recall that the indicator function of any rectangle in $\mathbb
R^d$ is a multiplier for all $p$, $1<p<\infty$ (see
Section~3,\,B); hence
$$
1_{[0, 2\pi]^d}\cdot(f\circ
h_1)\in\bigcap_{1<p<\infty}M_p(\mathbb R^d)
$$
for all $f\in K$. Since the homeomorphism $\varphi_1$ maps the
interval $[0, 2\pi]$ onto itself, we can regard it as a
self-homeomorphism of the circle $\mathbb T$. So, we can regard
$h_1$ as a self-homeomorphism of the torus $\mathbb T^d$.
Denote this self-homeomorphism of the torus by $h_2$. Clearly,
$f\circ h_2$ is a $2\pi$-periodic in each variable extension of
the function $1_{[0, 2\pi]^d}\cdot(f\circ h_1)$. It remains to
use the theorem on periodic extensions (see Section~3,\,D).
This completes the proof of Theorem 2.

\quad

\begin{center}
\textbf{5. Remarks and Open Problems}
\end{center}

1. As we stated in Introduction, the negative solution of
Luzin's problem implies a similar result for functions on the
real line, namely: there exists a bounded continuous
real-valued function $f$ on $\mathbb R$ such that $f\circ
h\notin B(\mathbb R)$ for every self-homeomorphism $h$ of
$\mathbb R$. One can easily verify this as follows. Let $f$ be
a function on $\mathbb T$ which provides the negative solution
of Luzin's problem. Regarding $f$ as a $2\pi$-periodic function
on $\mathbb R$, we can assume without loss of generality that
$f(0)=f(2\pi)=0$. Let $f_0$ be the function on $\mathbb R$ that
coincides with $f$ on the interval $[0, 2\pi]$ and vanishes
outside it. Suppose that $f_0\circ h\in B(\mathbb R)$ for some
self-homeomorphism $h$ of $\mathbb R$. Clearly, if a function
belongs to $B(\mathbb R)$, then so does every superposition of
this function with an affine mapping $l \colon \mathbb
R\rightarrow\mathbb R$; so, replacing $h$ by $h\circ l$ if
necessary, we can assume that $h$ maps $[0, 2\pi]$ onto itself,
and hence $f_0\circ h$ vanishes outside $[0, 2\pi]$. It follows
that $2\pi$ -periodic extension of $f_0\circ h$ to $\mathbb R$
is in $A(\mathbb T)$ (see Section 3, D), which is impossible,
since this extension has the form $f\circ h_1$, where $h_1$ is
a self-homeomorphism of the circle.

2. It is natural to consider an analogue of Luzin's problem in
the multidimensional case. Is it true that, given a real-valued
function $f\in C(\mathbb T^d), \,d\geq 2$, there exists a
self-homeomorphism $h$ of $\mathbb T^d$ such that $f\circ h\in
A(\mathbb T^d)$? Since the group of homeomorphisms of $\mathbb
T^d$ with $d\geq 2$ is more massive than that of $\mathbb T$
the question may have a positive answer for $d\geq 2$ despite
the fact that in the one-dimensional case it is answered in the
negative.

3. As in the introduction, let $G$ be either $\mathbb R^d$ or
$\mathbb Z^d$ and $\Gamma$ is $\mathbb R^d$ or $\mathbb T^d$,
correspondingly. Recall that a function $m\in L^\infty(\Gamma)$
is called a weak type $(1, 1)$ multiplier if the operator $Q$
defined by $\widehat{Qf}=m\widehat{f}, \,f\in L^1\cap L^2(G),$
is of weak type $(1, 1)$, i.e., satisfies the condition
$$
\mathrm{mes}_G\{t\in G : |Qf(t)|>\lambda\}\leq
c \|f\|_{L^1(G)}/\lambda,
\qquad \lambda>0,
$$
where $\mathrm{mes}_G$ is the Haar measure on $G$ (the Lebesgue
measure of a set in the case of $G=\mathbb R^d$ and the number
of elements in a set in the case of $G=\mathbb Z^d$). Let
$M_1^{weak}(\Gamma)$ denote the class of all such multipliers.
We have (using Marcinkiewicz' interpolation theorem)
$$
M_1(\Gamma)\subseteq M_1^{weak}(\Gamma)\subseteq \bigcap_{1<p<\infty}M_p(\Gamma).
$$
The authors do not know whether it is possible to improve
Theorems~1 and~2 so as to attain the condition $f\circ h\in
M_1^{weak}$ for all $f\in K$. The answer is unclear even for
the families that consist of one function: Given a bounded
real-valued function $f\in C(\mathbb R^d)$ or a real-valued
function $f\in C(\mathbb T^d)$ is there a homeomorphic change
of variable $h$ such that $f\circ h\in M_1^{weak}$? The
negative answer to this question in the one-dimensional case
would strengthen the result that solves Lusin's problem.

4. Let $d\geq 2$. Observe that the homeomorphism $h$ of
$\mathbb T^d$ in Theorem 2 can be chosen in the form
$$
h \colon (t_1, t_2,
\ldots, t_d)\rightarrow (\varphi(t_1), \varphi(t_2), \ldots,
\varphi(t_d)),
\eqno(15)
$$
where $\varphi$ is a self-homeomorphism of $\mathbb T$. As
concerns Theorem 1, it is clear from its proof that if we
impose stronger assumption on $f$ or on $K$, namely if we
assume that $f$ is bounded and uniformly continuous or,
respectively, that $K$ is a family of bounded functions which
is uniformly equicontinuous on the whole $\mathbb R^d, \,d\geq
2$, then the corresponding homeomorphism $h$ of $\mathbb R^d$
can be chosen in the form (15) with $\varphi$ being a
self-homeomorphism of $\mathbb R$. We do not know if the same
is true without the above stronger assumptions. Moreover, it is
unclear if every bounded real-valued function $f\in C(\mathbb
R^d)$ can be transformed into a multiplier by a homeomorphism
of the form $h \colon (t_1, t_2, \ldots, t_d)\rightarrow
(\varphi_1(t_1), \varphi_2(t_2), \ldots, \varphi_d(t_d))$,
where $\varphi_j$'s are self-homeomorphisms of $\mathbb R$
allowed to be different. It seems likely that the answer is
negative.

5. The well-known Beurling--Helson theorem [2] (see also [8,
Ch.~VI], [9]) states that if $\varphi$ is a continuous
self-mapping of the circle $\mathbb T$ satisfying the condition
$\|e^{in\varphi}\|_{A(\mathbb T)}=O(1)$, $n\in\mathbb Z$, then
$\varphi$ is linear (affine), i.e., $\varphi(t)=\nu
t+\varphi(0)$ where $\nu\in\mathbb Z$. The character of growth
of the norms $\|e^{in\varphi}\|_{A(\mathbb T)}$ for nontrivial
$\varphi$'s is unclear in many respects. Kahane conjectured
that the Beurling--Helson theorem can be considerably improved;
namely, he conjectured in [7] (see also [8, Ch.~VI], [9]) that
the conclusion of the Beurling--Helson theorem remains valid
even if the norms $\|e^{in\varphi}\|_{A(\mathbb T)}$ grow to
infinity but the growth is not very fast. He also conjectured
([7], [8]) that the condition $\|e^{in\varphi}\|_{A(\mathbb
T)}=o(\log|n|)$ already implies linearity of $\varphi$. The
first result in this direction was obtained in [17]; further
strengthening is obtained in [12], however the $o(\log |n|)$
-conjecture remains unproved.

Certain analogs of the Beurling--Helson theorem for the
algebras $M_p$ were obtained in [20] and [21]. Note that the
case $1<p<\infty$ differs from that of $p=1$, which corresponds
to the Wiener algebra; for example (see [5, Ch. I, Sec. 1.3]),
if $\varphi\colon \mathbb T\rightarrow \mathbb T$ is piecewise
linear, then for all $p, 1<p<\infty,$
$$
\|e^{in\varphi}\|_{M_p(\mathbb T)}=O(1), \qquad n\in \mathbb Z.
\eqno(16)
$$
To some extent the converse is also true [21, Theorem $2'$]: if
(16) holds for some $p\neq 2$, then there exists a closed set
$E(\varphi)\subset\mathbb T$ of Lebesgue measure zero such that
$\varphi$ is linear on the intervals complementary to
$E(\varphi)$, and the set of distinct slopes of $\varphi$ is
finite.

It is natural to ask how slow the norms
$\|e^{in\varphi}\|_{M_p}$ can grow in the case when $\varphi$
is nowhere (i.e., on no interval) linear. The answer to this
question can easily be extracted from results of this paper.
For $1<p<\infty$ the growth can be arbitrarily slow, namely:
\emph{given an arbitrary positive sequence $\gamma (n)$, $n=0,
1, 2, \ldots,$ with $\gamma (n)\rightarrow +\infty$, there
exists a nowhere linear self-homeomorphism $h$ of $\mathbb T$
such that
$$
\|e^{inh}\|_{M_p(\mathbb T)}=O(\gamma (|n|)), \qquad n\in \mathbb Z,
$$
for all $p$, $1<p<\infty$}.

To prove this assertion, we firstly observe that Theorems 1 and
2 can be supplemented by the estimate $\|f\circ h\|_{M_p}\leq
c(p, 1)\|f\|_{L^\infty}+c_p(K), \,f\in K$. In particular, if
$K$ is a compact set in the space $C(\mathbb T^d)$, then there
exists a homeomorphism $h \colon \mathbb T^d\rightarrow\mathbb
T^d$ such that
$$
\|f\circ h\|_{M_p(\mathbb T^d)}\leq c_K(p), \qquad 1<p<\infty, \quad f\in K.
\eqno(17)
$$
Indeed (see the end of the proof of Theorem 1), since
$g_\nu\overset{M_p}{\rightarrow} g$, it follows from (14) that
$\|g-g_1\|_{M_p}\leq 2 c_p(K)$, where $c_p(K)$ is the sum of
the series in~(9). At the same time $\|g_1\|_{M_p}\leq c(p,
1)\|g_1\|_{L^\infty}\leq c(p, 1)\|f\|_{L^\infty}$.

Secondly we observe that due to sufficient flexibility in the
construction of the net $\beta$ in the proof of Theorem 1, the
homeomorphism $h$ in the statement of Theorems 1 and 2 can be
made nowhere linear.

It remains to apply the above two observations to the family
$$
K=\bigg\{\frac{e^{int}}{\gamma(|n|)}, \quad n\in
\mathbb Z \bigg\}
$$
(see (17)).

\quad

\begin{center}
\textbf{References}
\end{center}

\flushleft
\begin{enumerate}

\item N. K. Bary, \emph{A treatise on trigonometric series},
    Vols. I, II, Pergamon Press, Oxford 1964.

\item A. Beurling, H. Helson, ``Fourier--Stieltjes transforms
    with bounded powers'', \emph{Math. Scand.,} 1 (1953),
    120--126.

\item R. E. Edwards and G. I. Gaudry, \emph{Littlewood-Paley
    and multiplier theory.} Springer-Verlag, Berlin--Heidelberg, 1977.

\item G. Goffman, T. Nishiura, D. Waterman,
    \emph{Homeomorphisms in Analysis}, Mathematical Surveys
    and Monographs v. 54, Amer. Math. Soc., 1997.

\item L. H\"ormander, ``Estimates for translation invariant
    operators in $L^p$ spaces'', \emph{Acta Math.}, 104
    (1960), 93--140.

\item M. Jodeit, ``Restrictions and extensions of Fourier
    multipliers'', \emph{Studia Math.}, 34 (1970), 215--226.

\item J.-P. Kahane, ``Transform\'ees de Fourier des fonctions
    sommables'', \emph{Proceedings of the Int. Congr. Math.,
    15-22 Aug., 1962, Stockholm, Sweden}, Inst.
    Mittag-Leffler, Djursholm, Sweden, 1963, pp.~114--131.

\item J.-P. Kahane, \emph{S\'erie de Fourier absolument
    convergentes}, Springer-Verlag, Berlin--Heidelberg--New
    York, 1970.

\item J.-P. Kahane, ``Quatre le\c cons sur les
    hom\'eomorphismes du circle et les s\'eries de Fourier'',
    in: \emph{Topics in Modern Harmonic Analysis, Vol. II,
    Ist. Naz. Alta Mat. Francesco Severi, Roma}, 1983,
    955--990.

\item  J.-P. Kahane, Y. Katznelson, ``Hom\'eomorphismes de
    cercle et s\'eries de Fourier absolument convergent'',
    \emph{Compt. Rend., Paris}, 292:4 (1981), S\'{e}rie I,
    271--273.

\item  J.-P. Kahane,  Y. Katznelson, ``S\'eries de Fourier
    des fonctions born\'ee'', \emph{Studies in Pure Math.},
    in Memory of Paul Tur\'an, Budapest, 1983, pp.~395--410.
    (Preprint, Orsay, 1978.)

\item S. V. Konyagin and I. D. Shkredov, ``A quantitative
    version of the Beurling--Helson theorem,'' \emph{Funct.
    Anal. Appl.} 49:2 (2015), 110--121.

\item G. Kozma, ``Random homeomorphisms and Fourier
    expansions
    --- the pointwise behaviour'', \emph{Israel Journal of
    Mathematics}, 139 (2004), 189--213.

\item G. Kozma and A. M. Olevskii, ``Random homeomorphisms
    and Fourier expansions'', \emph{Geometric and Fnnctional
    Analysis}, 8 (1998), 1016--1042.

\item V. V. Lebedev, ``Change of variable and the rapidity of
    decrease of Fourier coefficients'',
    \emph{Matematicheski\v{\i} Sbornik}, 181:8 (1990),
    1099--1113 (in Russian). English transl.:
    \emph{Mathematics of the USSR-Sbornik}, 70:2 (1991),
    541--555. English transl. corrected by the author is
    available at: https://arxiv.org/abs/1508.06673

\item V. V. Lebedev, ``Torus homeomorphisms, Fourier
    coefficients, and integral smoothness'', \emph{Russian
    Mathematics (Iz. VUZ)}, 36:12 (1992), 36--41.

\item V. V. Lebedev, ``Absolutely convergent Fourier series.
    An improvement of the Beurling--Helson theorem,''
    \emph{Funct. Anal. Appl.} 46:2 (2012), 121--132.

\item V. Lebedev, ``The Bohr--P\'al theorem and the Sobolev
    space $W_2^{1/2}$\,'', \emph{Studia Mathematica}, 231:1
    (2015), 73--81.

\item  V. V. Lebedev, ``A short and simple proof of the
    Jurkat--Waterman theorem on conjugate functions'',
    \emph{Functional Analysis and Its Applications}, 51:2
    (2017), 148--151.

\item V. Lebedev and A. Olevski\v{\i}, ``$C^1$ changes of
    variable: Beurling-Helson type theorem and H\"ormander
    conjecture on Fourier multipliers'', \emph{Geometric and
    Functional Analysis (GAFA)} 4 (1994), 213--235.

\item V. V. Lebedev, A. M. Olevski\v{\i}, ``$L^p$ -Fourier
    multipliers with bounded powers'', \emph{Izvestya:
    Mathematics}, 70:3 (2006), 549--585.

\item A. M. Olevski\v{\i}, ``Change of variable and absolute
    convergence of Fourier series'', \emph{Soviet Math.
    Dokl.}, 23 (1981), 76--79.

\item A. M. Olevski\v{\i}, ``Modifications of functions and
    Fourier series'', \emph{Russian Math. Surveys}, 40
    (1985), 181--224.

\item G. T. Oniani, ``Topological characterization of a set
    of continuous functions whose conjugate functions are
    continuous and have bounded variation'', \emph{Trudy
    Tbiliss. Mat. Inst. Razmadze Akad. Nauk Gruzin. SSR}, 86
    (1987), 110--113 (in Russian). English transl.:
    \emph{Integral Operators and Boundary Properties of
    Functions. Fourier Series}, 155--161, Nova Science, New
    York, 1992.

\item A. A. Saakjan, ``Integral moduli of smoothness and the
    Fourier coefficients of the composition of functions'',
    \emph{Matematicheski\v{\i} Sbornik}, 110(152):4(12)
    (1979), 597--608 (in Russian). Engl. transl.:
    \emph{Mathematics of the USSR-Sbornik}, 38 (1981),
    549--561.

\item P. Sj\"{o}gren and P. Sj\"{o}lin, ``Littlewood-Paley
    decompositions and Fourier multipliers with singularities
    on certain sets.'', \emph{Ann. Inst. Fourier, Grenoble},
    31 (1981), 157--175.

\end{enumerate}

\quad

\qquad Vladimir Lebedev\\
\qquad School of Applied Mathematics\\
\qquad National Research University Higher School of Economics\\
\qquad 34 Tallinskaya St.\\
\qquad Moscow, 123458 Russia\\
\qquad E-mail address: \emph{lebedevhome@gmail.com}

\quad

\qquad Alexander Olevskii\\
\qquad School of Mathematical Sciences\\
\qquad Sackler Faculty of Exact Sciences\\
\qquad Tel Aviv University\\
\qquad Tel Aviv, 69978 Israel\\
\qquad E-mail address: \emph {olevskii@yahoo.com}

\end{document}